\documentclass[a4paper]{eccomas_2022}
\usepackage{graphicx}
\usepackage{amsmath}
\usepackage{amsfonts}
\usepackage{amssymb}

\usepackage{graphicx}
\usepackage{amsmath}
\usepackage{amsfonts}
\usepackage{amssymb}

\newcommand{\intd}{\textrm{\;d}}
\newcommand{\ddt}{\frac{\textrm{d}}{\textrm{d} t}}

\DeclareMathOperator{\rank}{rank}

\newcommand{\h}{\hspace{1 em}}

\newcommand{\vv}[1]{\ensuremath{\mathbf{#1}}} 
\usepackage{tikz}
\usetikzlibrary{shapes.geometric}
\usetikzlibrary{positioning}
\usetikzlibrary{positioning,shapes.multipart, fit,backgrounds,calc}

\usepackage[utf8]{inputenc}
\usepackage[T1]{fontenc}
\usepackage{ae}
\usepackage[english]{babel}
\usepackage{amsmath,amssymb,amstext,amsfonts}
\usepackage{units}
\usepackage[squaren]{SIunits}
\usepackage{bm}
\usepackage{graphicx}
\usepackage{pgf}
\usepackage{pgfplots}

\usepackage{tikz}
\usetikzlibrary{arrows}
\usetikzlibrary{shapes}
\usetikzlibrary{decorations.text}
\usepackage{xcolor}
\usepackage{color}
\usepackage{tikzscale}
\usepackage{pgfplots}

\usetikzlibrary{shapes.geometric}
\usetikzlibrary{positioning}
\usetikzlibrary{positioning,shapes.multipart, fit,backgrounds,calc}

\usepackage{nicefrac}
\usepackage{array,tabularx}
\usepackage{multirow}
\newcolumntype{L}[1]{>{\raggedright\arraybackslash}p{#1}}
\newcolumntype{C}[1]{>{\centering\arraybackslash}p{#1}}
\newcolumntype{R}[1]{>{\raggedleft\arraybackslash}p{#1}}
\usepackage{booktabs}
\usepackage{mathtools}

\usepackage{geometry}
\usepackage{pdfpages}

\usepackage[colorinlistoftodos,prependcaption,textsize=tiny]{todonotes}
\usepackage{caption}
\usepackage{subcaption}
\usepackage{framed}
\usepackage[pdfencoding=auto]{hyperref}
\usepackage{placeins}

\usepackage[amsthm]{ntheorem}
\theoremheaderfont{\normalfont\bfseries}

\theoremstyle{remark}

\theoremstyle{break}

\newcommand{\bitem}{\begin{itemize}}
\newcommand{\eitem}{\end{itemize}}
\newcommand{\bsitem}{\begin{itemize*}}
\newcommand{\esitem}{\end{itemize*}}
\newcommand{\benum}{\begin{enumerate}}
\newcommand{\eenum}{\end{enumerate}}

\newcommand{\bdesc}{\begin{description}}
\newcommand{\edesc}{\end{description}}

\newcommand{\mcM}{\mathcal{M}}

\newcommand{\mbR}{\mathbb{R}} 

\newcommand{\beq}{\begin{eqnarray*}}
\newcommand{\eeq}{\end{eqnarray*}}
\newcommand{\beqn}{\begin{eqnarray}}
\newcommand{\eeqn}{\end{eqnarray}}

\newcommand{\bpmat}{\begin{pmatrix}} 
\newcommand{\epmat}{\end{pmatrix}}
\newcommand{\bbmat}{\begin{bmatrix}} 
\newcommand{\ebmat}{\end{bmatrix}}

\ifnum 0=1

\def\XXint#1#2#3{{\setbox0=\hbox{$#1{#2#3}{\int}$ }
\vcenter{\hbox{$#2#3$ }}\kern-.6\wd0}}

\fi
\theoremstyle{plain}

\theoremstyle{plain}

\title{Momentum-conserving ROMs for the incompressible Navier-Stokes equations}

\author{Henrik K. E. Rosenberger$^1$ and Benjamin Sanderse$^2$}

\heading{Henrik K. E. Rosenberger and Benjamin Sanderse}

\address{$^{1}$ Centrum Wiskunde \& Informatica \\ 
Science Park 123, Amsterdam, The Netherlands \\
email: henrik.rosenberger@cwi.nl, www.cwi.nl/people/henrik-rosenberger
\and
$^{2}$  Centrum Wiskunde \& Informatica \\
Science Park 123, Amsterdam, The Netherlands \\
email: b.sanderse@cwi.nl, www.thinkingslow.nl}


\abstract{

Projection-based model order reduction of an ordinary differential equation (ODE) 
results in a projected ODE. 
Based on this ODE, an existing reduced-order model (ROM) for finite volume discretizations satisfies the underlying conservation law over arbitrarily chosen subdomains.
However, 
this ROM
does not satisfy the projected ODE exactly but introduces an additional perturbation term. 
In this work, we propose a novel ROM with the same subdomain conservation properties which indeed satisfies the projected ODE exactly.

We apply this ROM to the incompressible Navier-Stokes equations and show with regard to the mass equation how the novel ROM can be constructed to satisfy algebraic constraints.

Furthermore, we show that the resulting mass-conserving ROM allows us to derive kinetic energy conservation and consequently nonlinear stability, which was not possible for the existing ROM due to the presence of the perturbation term. 

 }

\begin{document}

\section{INTRODUCTION}
For many applications such as optimization, uncertainty quantification and real-time simulation of large-scale systems, classical numerical methods such as finite elements or finite volumes are prohibitively expensive. To mitigate the computational costs, several techniques have been developed to reduce the complexity of the full order models (FOM) developed via such classical methods, yielding so-called reduced order models (ROM). Many of these techniques are projection-based, i.e., the FOM is projected onto a lower-dimensional subspace spanned by a reduced basis. Such a reduced basis can be obtained, e.g., via Proper Orthogonal Decomposition of FOM snapshots \cite{volkwein2011model}.

While ROMs have been developed successfully for elliptic and parabolic systems of partial differential equations, deriving stable ROMs for hyperbolic systems is challenging \cite{afkham2020conservative}. Stability problems of these ROMs are often ascribed to lacking structure preservation. While many FOMs feature physical structures of the underlying systems, many model reduction techniques do not preserve such structures \cite{afkham2020conservative}. 

Therefore, many approaches to obtain stable reduced order models refine model reduction techniques such that they do preserve such physical structures \cite{hesthaven2022reduced}. Peng et al.\ \cite{peng2016symplectic} propose a POD-Galerkin ROM for Hamiltonian systems that preserve the systems' symplecticity. Chan \cite{chan2020entropy} constructs ROMs for nonlinear conservation laws that exhibit a semi-discrete entropy dissipation. Sanderse \cite{sanderse2020non} achieves nonlinear stability of a POD-Galerkin ROM for the incompressible Navier-Stokes equations by preserving the kinetic energy evolution and additionally conserves the global momentum.


Carlberg et al.\ \cite{carlberg2018conservative}
address another aspect of structure preservation specific to finite volume discretizations. By construction, these discretizations satisfy the underlying conservation law over all finite volumes. Projection-based ROMs, on the other hand, satisfy this conservation law generally over none of the finite volumes. While it is in most cases impossible to define ROMs that satisfy the conservation law over all finite volumes, Carlberg et al.\ propose ROMs that satisfy the conservation law at least over a few subdomains. These ROMs are formulated via a constrained optimization problem (COP) and are shown to be more accurate than ROMs without subdomain conservation.
Furthermore, the COP framework has been extended by
Schein et al.\ \cite{schein2021preserving} to include other types of constraints.

In this work, we combine the concept of conservation over subdomain in \cite{carlberg2018conservative} and the nonlinear stable model reduction in \cite{sanderse2020non}. Unfortunately, the subdomain conservative ROM in \cite{carlberg2018conservative} is not suitable for this purpose. As shown in \cite{carlberg2018conservative}, the solution to the COP can equivalently be expressed by an ODE. This ODE, however, differs from the conventional Galerkin ROM ODE by a perturbation term. This perturbation term impedes the derivation of nonlinear stability as described in \cite{sanderse2020non}. Therefore, we propose a novel subdomain conservative ROM which does not introduce any perturbation terms.

We apply the novel ROM to the energy-conserving finite volume discretization of the incompressible Navier-Stokes equations proposed in \cite{sanderse2020non}. Regarding the momentum equation, we can directly apply the novel ROM. Regarding the mass equation, however, we have to extend the novel ROM to treat algebraic constraints. The resulting mass-conserving ROM is shown to mimick the kinetic energy evolution of the FOM. This evolution states that the global kinetic energy does not increase over time, hence implying nonlinear stability of the ROM.



This article is structured as follows.
First, we summarize the subdomain conservative ROM 
proposed in
\cite{carlberg2018conservative} 
in Section \ref{sec:carlberg cons for general cons law}.
Then, we derive our novel subdomain conservative ROM in Section \ref{sec:novel approach}, in the context of a general conservation law. In Section \ref{sec:mom cons over subdoms for ins}, we apply our ROM to the energy-conserving finite volume discretization of the incompressible Navier-Stokes equations proposed in \cite{sanderse2020non} and show that the ROM preserves the kinetic energy evolution of the FOM.
We summarize our findings in Section \ref{sec:conclusion} and give an outlook on possible uses of the proposed ROMs in Section \ref{sec:outlook}.

\section{SUBDOMAIN CONSERVATION FOR GENERAL CONSERVATION LAWS}
\label{sec:carlberg cons for general cons law}

In this section, we summarize the concept of subdomain conservation and the subdomain conservative ROM described in \cite{carlberg2018conservative}. For convenience only, we use a simplified notation and consider only scalar-valued conservation laws.

\subsection{Subdomain conservation for finite volume discretizations}
In the understanding of Carlberg et al.\ \cite{carlberg2018conservative},
a 
scalar-valued
quantity $u$ is 
conserved over a domain $\omega$ if it satisfies the integral form of a 
conservation law over this domain
\begin{align}
	\ddt \int_\omega u \intd V + \int_{\partial \omega} g\cdot 
	n\intd S = \int_\omega s \intd V .
	\label{eq:integral conservation law} 
\end{align}
where $g=g(u)$ describes the flux across the boundary and $s$ is a term 
describing sources and 
sinks.

For example,
a quantity $u$ satisfying the differential form of a conservation 
law
\begin{align}
	\frac{\partial}{\partial t}u + \nabla \cdot g = s \hspace{2em} \textrm{ in 
	} 
	\Omega,
	\label{eq:conservation law}
\end{align}
is conservative over any subdomain of $\Omega$, because \eqref{eq:conservation 
law} 
 implies \eqref{eq:integral conservation law} for all 
 $\omega\subset\Omega$.
	
Let us now consider a finite volume discretization of \eqref{eq:conservation 
law} on a mesh $\mcM=\{\Omega_j\}_{j=1}^{N_\mcM}$ 
that completely covers $\Omega$.
Then, we can define the time-dependent 
state-vector $\vv{u_h} = \begin{bmatrix}	u_1	\;\dots\;  u_j \; \dots\; u_{N_{\mcM}}
\end{bmatrix}^T$ of cell averages
\begin{align}
	u_{j} = \frac{1}{\left|\Omega_j\right|}\int_{\Omega_j} u \intd V 
	\hspace{1em} 
	j=1,\dots,N_\mcM ,
	\label{eq:fv cell averages}
\end{align}
and the time-dependent vector $\vv f = \begin{bmatrix}
	-f^g_1+f^s_1, \;	\dots \; -f^g_j+f^s_j, \; \dots \; 
	-f^g_{N_{\mcM}}+f^s_{N_{\mcM}}
\end{bmatrix}^T$ 
consisting of the net flux 
	\begin{align}
	f^g_j &= \frac{1}{\left|\Omega_j\right|} \int_{\partial\Omega_j} g^{FV} 
	\cdot 
	n \intd S ,
	\label{eq:fv cell net flux} \\
\intertext{and the net sinks and sources}
	f^s_j &= \frac{1}{\left|\Omega_j\right|} \int_{\Omega_j} s^{FV} \intd V,
	\label{eq:fv cell net sinks and sources}
\end{align}
of each cell, divided by the respective cell-volume. Here, $g^{FV}$ and 
$s^{FV}$ are approximations of $g$ and $s$ in \eqref{eq:conservation law}.  

Together, the state vector $\vv{u_h}(t)$ and the RHS vector $\vv f(\vv{u_h},t)$ form the 
ordinary 
differential equation (ODE) system
\begin{align}
	\ddt \vv{u_h} =\vv f(\vv{u_h},t) ,
	\label{eq:general cons law fom ode}
\end{align}
which we can interpret as a system of conservation laws in integral form over 
all finite volume cells. 

Even more, the vector $\vv{u_h}$ represents a function on $\Omega$ which is constant over each finite volume.
This function
is subdomain conservative over all finite volume cells 
with respect to
the 
approximated fluxes $g^{FV}$, and sink and source terms $s^{FV}$. Therefore, we 
also denote $\vv{u_h}$ itself as subdomain conservative over all finite volume cells.

In the following, it will be useful to express the ODE \eqref{eq:general cons law fom ode}
equivalently as 
\begin{align} 
	\vv r\left(\ddt \vv{u_h}, \vv{u_h},t\right) = 0,
\end{align}
with the residual $\vv r(\vv v,\vv w,t) = \vv v - \vv f(\vv w,t)$.

\subsection{Subdomain conservation for Galerkin projection ROMs}
Let us now investigate the subdomain conservation properties of a 
model 
ROM approximating the FOM
\eqref{eq:general cons law fom ode}.
Given an orthonormal basis $\Phi\in\mbR^{N_\mcM\times p}, \;p\ll N_\mcM$, we 
consider the ROM 
\begin{align}
	 \vv{u_r}(t) = \Phi \vv a(t)  \approx \h \vv{u_h}(t),
	\label{eq:general rom} 
\end{align}
that satisfies the Galerkin projection of the FOM ODE 
\eqref{eq:general cons law fom ode}
\begin{align}
	\Phi^T\frac{\intd}{\intd t}\vv{u_r} = \frac{\intd}{\intd t}\vv a = 
	\Phi^T\vv f(\Phi\vv a,t) \hspace{1em} \Leftrightarrow \hspace{1em} 
	\Phi^T \vv r\left(\Phi \ddt \vv a, \Phi \vv a,t\right)=0.
	\label{eq:general rom ode}
\end{align}
As shown in \cite{carlberg2017galerkin}, this 
ODE system, equipped with some initial condition $\vv a(0)= \vv a_0$, can be 
equivalently described by the optimization problem
\begin{align}
	\ddt \vv a= \arg \min_{\vv b \in \mbR^p}\|r\left(\Phi\vv b,\Phi\vv a, 
	t\right)\|_2,
	\label{eq:optimi}
\end{align}
and the same initial condition.
Hence, the ROM $\vv{u_r}(t)$ minimizes the violation of the cell-wise conservation 
laws\,\eqref{eq:general cons law fom ode} but this minimum is not guaranteed to 
vanish.
Consequently, the ROM $\vv{u_r}(t)$ is not guaranteed to be subdomain 
conservative over any finite volume.

Hence, instead of aiming at subdomain conservation over all finite volumes, 
Carlberg et al.\
propose to require conservation over a set of a few subdomains each 
consisting 
of one or more finite volumes. 

\begin{figure}
	\centering
	\scalebox{0.8}{	
%
%
\begin{tikzpicture}[thick,black]
\tikzset{edot/.style={fill=black,diamond}}

\foreach \y in {0,1,...,6}
{
\draw (0,\y) -- (6,\y);
\draw (\y,0) -- (\y,6);
}

\foreach \y in {1,3}
{
	\draw [red,line width=2mm] (1,\y) -- (3,\y);
	\draw [red,line width=2mm] (\y,1) -- (\y,3);
}

\draw [blue,line width=2mm] (2,6)--(5,6);
\draw [blue,line width=2mm] (5,6)--(5,5);
\draw [blue,line width=2mm] (5,5)--(6,5);
\draw [blue,line width=2mm] (6,5)--(6,4);
\draw [blue,line width=2mm] (6,4)--(4,4);
\draw [blue,line width=2mm] (4,4)--(4,5);
\draw [blue,line width=2mm] (4,5)--(2,5);
\draw [blue,line width=2mm] (2,5)--(2,6);

\foreach \x in {4}
{	\foreach \y in {1}
	{
		\draw [green, line width=1mm] (\x,\y)--(\x,\y+1);
		\draw [green, line width=1mm] (\x+1,\y)--(\x+1,\y+1);
		\draw [green, line width=1mm] (\x,\y)--(\x+1,\y);
		\draw [green, line width=1mm] (\x,\y+1)--(\x+1,\y+1);
	}
}

\foreach \x in {2}
{	\foreach \y in {2}
	{
		\draw [green, line width=1mm] (\x,\y)--(\x,\y+1);
		\draw [green, line width=1mm] (\x+1,\y)--(\x+1,\y+1);
		\draw [green, line width=1mm] (\x,\y)--(\x+1,\y);
		\draw [green, line width=1mm] (\x,\y+1)--(\x+1,\y+1);
	}
}

\foreach \x in {0}
{	\foreach \y in {4}
	{
		\draw [green, line width=1mm] (\x,\y)--(\x,\y+1);
		\draw [green, line width=1mm] (\x+1,\y)--(\x+1,\y+1);
		\draw [green, line width=1mm] (\x,\y)--(\x+1,\y);
		\draw [green, line width=1mm] (\x,\y+1)--(\x+1,\y+1);
	}
}

\end{tikzpicture}
 }
	\caption{Example of three subdomains.}
	\label{fig: subdoms}
\end{figure}
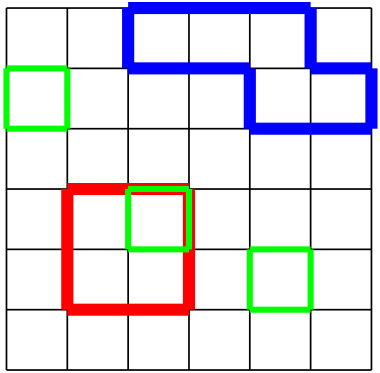

To this 
end, we consider
a decomposed mesh $\bar\mcM$ of $N_{\bar\mcM}$ arbitrary subdomains $\bar 
\Omega_k = 
\cup_{j\in 
	S_k} \Omega_j$ of  
finite volumes $\Omega_j$ where $S_k\subset \{1,\dots,N_\mcM\}$. These subdomains $\bar\Omega_h$ can overlap and do not need to be connected. An example of such a decomposed mesh is depicted in 
Fig.\ \ref{fig: subdoms}. 
On this grid, 
we define, similarly to the finite volume quantities \eqref{eq:fv cell averages}, \eqref{eq:fv cell net flux} and \eqref{eq:fv cell net sinks and sources}, a time-dependent state 
vector $\vv{\bar u_h} = \begin{bmatrix}
	\bar u_1	\;\dots\; \bar u_k \; \dots\; \bar u_{N_{\bar \mcM}}
\end{bmatrix}^T$ of subdomain averages
\begin{align}
	\bar u_k = \frac{1}{\left|\bar\Omega_k\right|}\int_{\bar\Omega_k} u \intd V 
	= 
	\frac{1}{\left|\bar \Omega_k\right|}\sum_{j\in S_k} \int_{\Omega_j} u \intd 
	V = 
	\frac{1}{\left|\bar \Omega_k\right|}\sum_{j\in S_k} 
	\left|\Omega_j\right|u_j ,
\end{align}
and analogously a time-dependent RHS vector
$\vv{\bar f}
,$
consisting of the net flux 
and the net sinks and sources
of each subdomain, divided by the respective subdomain volume. 
%
We
can write $\vv{\bar u_h}$ and $\vv{\bar f}$ in matrix vector notation as
\begin{align}
	\vv{\bar u_h} &= C^T \vv{u_h} & \vv{\bar f} &= C^T \vv f,
\end{align}
for some $C\in\mbR^{N_{\mcM}\times N_{\bar\mcM}}$.

As a result, subdomain conservation over all subdomains of $\bar\mcM$ is
described by the 
$C^T$-premultiplied ODE system
\begin{align}
	C^T \frac{d}{dt} \vv{u_h} = C^T\, \vv f \hspace{1em} \Leftrightarrow \hspace{1em} 
	C^T\,\vv r\left(\frac{d}{dt}\vv{u_h},\vv{u_h},t\right)=0.
	\label{eq:conservation over subdomains}
\end{align}

Based on the Galerkin projection ROM \eqref{eq:general rom ode}, we want to construct a ROM that also satisfies the subdomain conservation constraint \eqref{eq:conservation over subdomains}. 
In the following section, we summarize the subdomain conservative ROM proposed in \cite{carlberg2018conservative}. Our novel subdomain conservative ROM is presented in Section \ref{sec:novel approach}.

\subsection{Existing approach for subdomain conservative ROMs: constrained optimization problem formulation}
Carlberg et al.\ \cite{carlberg2018conservative} suggest 
to use the interpretation of 
the ROM ODE \eqref{eq:general rom ode} as minimization problem \eqref{eq:optimi}. Adding the conservation constraint \eqref{eq:conservation over subdomains} to this optimization problem, we get the constrained optimization problem
\begin{align}
    &\text{minimize}_{\,\vv b\in\mbR^{p}} & \|r(\Phi \vv b, \Phi \vv a, t)&\|_2 
    \label{eq:constr minimi minimi}
    \\
    &\text{subject to} \h & C^T\vv r(\Phi \vv b,\Phi \vv a, t) &= 0.
    \label{eq:constr minimi constr}
\end{align}

A sufficient condition for feasibility of this COP is given in \cite[Proposition 5.1]{carlberg2018conservative}. If the COP is feasible, the solution is equivalently expressed by the ODE
\begin{align}
    \ddt \vv a = \Phi^T \vv f(\Phi \vv a,t) + \vv{f^*}(\Phi \vv a,t)
\end{align}
with the perturbation term
\begin{align}
    \vv f^*(\Phi \vv a,t) = (C^T\Phi)^+[C^T-C^T\Phi\Phi^T]\vv f(\Phi \vv a,t),
    \label{eq:perturbation term}
\end{align}
and the Moore-Penrose inverse $(C^T\Phi)^+$ of $C^T\Phi$. 

Except for the perturbation term, this ODE is equivalent to the conventional Galerkin ROM ODE \eqref{eq:general rom ode}. The perturbation term, however, lacks a clear physical motivation and impedes the preservation of the kinetic energy evolution as we will see in Section \ref{sec:kinetiv energy}. Therefore, we propose a subdomain conservative ROM without any perturbation terms in the next section.



\section{NOVEL APPROACH FOR SUBDOMAIN CONSERVATIVE ROMS: BASIS MODIFICATION}
\label{sec:novel approach}

\subsection{Novel approach}
\label{sec:novel rom}
The basic idea of our novel approach is
to merge the ROM ODE \eqref{eq:general rom ode} and the conservation constraint \eqref{eq:conservation over subdomains} in one linear system for $\ddt \vv{ u_r}$,
\begin{align}
    \left[\Phi \;\; C\right]^T \ddt \vv{u_r} = \left[\Phi \;\; C\right]^T \vv f(\vv{u_r},t).
    \label{eq:merged system}
\end{align}
To satisfy this linear system, we define an orthogonal basis $\tilde \Phi\, \in \mbR^{N_\mcM\times H}, \; H := \rank([\Phi \;\; C])$ that spans the same linear subspace as $[\Phi \;\; C]$ and equip our novel ROM with this basis,
\begin{align}
    \vv{u_r}(t) = \tilde \Phi \vv{\tilde a}(t), \; \vv{\tilde a}(t) \in\mbR^{H} 
    \label{eq:u r expansion tilde phi}.
\end{align}
Then, the Galerkin ROM ODE
\begin{align}
    \tilde \Phi^T \ddt \tilde \Phi \vv{\tilde a} = \ddt \vv{\tilde a} = \tilde \Phi^T f(\tilde \Phi \vv{\tilde a},t),
    \label{eq:merge linear system2}
\end{align}
satisfies the
system \eqref{eq:merged system}. Hence, this ROM is subdomain conservative without introducing any perturbation terms.



\subsection{View on POD bases}

A popular method to construct ROM bases is POD, which exhibits an optimality property:
given a snapshot matrix $X=[\vv{u_h^0} \;\; \dots \;\; \vv{u_h^{K-1}}]\in\mbR^{N_\mcM \times K}$, a POD basis $\Phi_\mathrm{POD}\in\mbR^{N_\mcM\times p}$ is known to minimize the sum of best approximation errors of the snapshots in $X$ among all orthogonal matrices in $\mbR^{N_\mcM\times p}$ \cite{volkwein2011model}. 

On the other hand, the key element of our novel approach is the modification of the ROM basis in consideration of the predefined subdomains.
In order that the Galerkin ROM ODE \eqref{eq:merge linear system2} includes the subdomain conservation constraint \eqref{eq:conservation over subdomains}, the span of the basis $\tilde \Phi$ must include the span of $C$. This requirement can be expressed as the equivalence of $C$ and the orthogonal projection of $C$ onto $\tilde \Phi$, $C = \tilde\Phi\tilde\Phi^T C$. Combining the minimization idea of POD and the subdomain conservation constraint, we find the constrained optimization problem for the ROM basis $\tilde \Phi$
\begin{align}
    &\text{minimize}_{\;\Xi\in\mbR^{N_\mcM\times q}}  &\sum_{j=0}^{K-1} \| \vv{u_h^j} - \Xi\Xi^T \vv{u_h^j} &\|^2 \\
    &\text{subject to} & \Xi^T\Xi &= I \\
    &\text{and} & C &= \Xi\Xi^TC.
\end{align}
The orthogonal basis that spans the same subspace as 
$[\Phi_\mathrm{POD} \;\; C]$ 
does generally not solve this constrained optimization problem.

A solution to this constrained optimization problem is proposed in \cite{xiao2013constrained}.
This solution is defined as
\begin{align}
    \tilde \Phi = [Q_1 \;\; Q_2 V],
    \label{eq:def tilde phi}
\end{align}
where $Q_1\in\mbR^{N_\mcM\times H_C}$ and $Q_2\in\mbR^{N_\mcM\times (N_\mcM - H_C})$ with $H_C = \rank(C)$ are obtained from the QR decomposition of $C^T$,
\begin{align}
    C = [Q_1\;\;Q_2]\begin{bmatrix}
    R_1 \\ 0
    \end{bmatrix}, \; R_1\in\mbR^{H_C\times N_{\bar \mcM}}
\end{align}
and $V\in\mbR^{(N_\mcM-H_C)\times(q-H_C)}$ consists of the first $q-H_C$ modes of the POD basis of $Q_2^TX$. 

In view of Section \ref{sec:novel rom}, $V$ can be interpreted as the initial basis $\Phi$.
Observe that in contrast to $\Phi_\mathrm{POD}$, the matrix $V$ depends on $C$. 




\subsection{View on least-squeares Petrov-Galerkin ROMs}
In \cite{carlberg2018conservative}, Carlberg et al.\ discuss besides the Galerkin ROM also a least-squares Petrov-Galerkin (LSPG) ROM. 
The two ROMs differ in the order of performing the ROM approximation and the time discretization in their derivations.
While the Galerkin ROM is obtained by introducing the ROM approximation on the time-continuous model and discretizing in time afterwards, the LSPG ROM is obtained by first discretizing in time and then performing the ROM approximation on the time-discrete model.

So far, we have introduced our novel approach only for the Galerkin ROM. In fact, our approach does not directly work in the context of the LSPG ROM. For the LSPG ROM trial and test basis are not the same. As a consequence, the choice of the basis $\eqref{eq:def tilde phi}$ as trial basis does generally not imply
that the span of the test basis includes the span of $C$.
Hence, the subdomain conservation constraint \eqref{eq:conservation over subdomains} is generally not satisfied.


\subsection{Subdomain conservation for invariants}


We want to highlight a property of our novel ROM for the special case of invariants.
Invariants are quantities that are constant over time, e.g., the global integral of a conservation quantity with periodic boundary conditions and in absence of external forces.
If constraints in the constraint matrix $C$ describe such invariants,
then the coefficients in $\vv{\tilde a}(t)$ in the expansion \eqref{eq:u r expansion tilde phi} that correspond to those invariant constraints are constant. Hence, we do not need to integrate these coefficients via the ODE \eqref{eq:merge linear system2}, but can compute their values from the initial conditions. This observation can be seen as a physical justification for a ROM Ansatz with a constant offset, $\vv{u_r}(t) = \tilde \Phi \vv{\tilde a}(t) + \vv{u_r^0}$.


\section{MOMENTUM CONSERVATION OVER SUBDOMAINS FOR THE INCOMPRESSIBLE NAVIER-STOKES EQUATIONS}
\label{sec:mom cons over subdoms for ins}


\subsection{Introduction to the FOM}
We consider the energy-conserving finite volume discretization of the incompressible Navier-Stokes equations with periodic boundary conditions proposed in \cite{sanderse2020non},
\begin{align}
    M_h \vv{V_h}(t) &= 0, 
    \label{eq:mass eq} \\
    \Omega_h\ddt \vv{V_h}(t) &= \vv{F_h^{CD}}(\vv{V_h}(t),t) - G_h\vv{p_h}(t).
    \label{eq:mom eq}
\end{align}
The vectors $\vv{V_h}(t)\in\mbR^{N_V}$ and $\vv{p_h}\in\mbR^{N_p}$ describe the velocity and the pressure, respectively. The matrices $M_h\in\mbR^{N_p\times N_V}$ and $G_h\in\mbR^{N_V\times N_p}$ are the discretizations of the divergence and the gradient operator, respectively, and satisfy the duality property $M_h = -G_h^T$.
The matrix $\Omega_h\in\mbR^{N_V\times N_V}$ is diagonal with the sizes of the finite volumes on its diagonal, so symmetric positive definite. The term $\vv{F_h^{CD}}(\vv{V_h}(t),t)\in\mbR^{N_V}$ comprises convective, diffusive and body force contributions.

While the momentum equation \eqref{eq:mom eq} is an ODE, the mass equation \eqref{eq:mass eq} is an algebraic equation. The subdomain conservation ROM method described in Section \ref{sec:novel approach}, however, addresses conservation laws of ODE structure only. Therefore, we apply the subdomain conservation method only to the momentum equation and deal with the mass equation separately in Section \ref{sec:mass cons for arb velo basis}.

\subsection{Momentum conservation over subdomains}
To construct a ROM that conserves momentum over a set of predefined subdomains represented by the matrix $C\in\mbR^{N_V\times N_{\hat \mcM}}$,
we define the approximation of the velocity 
\begin{align}
    \vv{V_r}(t) = \phi \vv a(t) \h \approx \vv{V_r}(t), \h \phi\in\mbR^{N_V\times R_V}, \h \vv a(t)\in\mbR^{R_V}, \h R_V\ll N_V.
    \label{eq:velo expansion}
\end{align}
In contrast to Section \ref{sec:novel approach}, where we introduced the orthogonal basis $\tilde \Phi$ with respect to the $L^2$-inner product, we construct the velocity basis $\phi$ to be orthogonal with respect to the $\Omega_h$-weighted inner product, i.e.\ $\phi^T\Omega_h\phi = I$. As a result, the definition of $\phi$,
\begin{align}
    \phi = [\tilde Q_1\;\; W],
    \label{eq:def phi}
\end{align}
differs slightly from the definition of $\tilde \Phi$ \eqref{eq:def tilde phi}.
Here, $\tilde Q_1 = \Omega_h^{-1/2}Q_1$ with the QR decomposition of $\Omega_h^{-1/2}C$,
\begin{align}
    \Omega_h^{-1/2}C = [Q_1\;\;Q_2]\begin{bmatrix}
    R_1 \\ 0
    \end{bmatrix} = Q_1R_1,
\end{align}
and $W$ consists of the first $R_V-\rank(C)$ columns of the $\Omega_h$-orthogonal POD basis $\tilde U$ of a snapshot matrix $X = [\vv{V_h^0} \;\; \dots \;\; \vv{V_h^{K-1}}]\in\mbR^{N_V\times K}$ that is $\Omega_h$-orthogonal to $\tilde Q_1$. In detail: 
\begin{align}
    &\text{We project the snapshot matrix $X$} & & \nonumber \\
    &\text{onto the subspace $\Omega_h$-orthogonal to $\tilde Q_1$:}
    & \tilde X &= [I-\tilde Q_1\tilde Q_1^T\Omega_h]X. \\
    &\text{We transform to include the $\Omega_h$-inner product:} & \hat X &= \Omega_h^{1/2}\tilde X. \\
    &\text{We compute the SVD of $\hat X$:} & \hat X &= \hat U\hat\Sigma \hat V^T. \\
    &\text{We transform back:} & \tilde U &= \Omega_h^{-1/2}\hat U. \\
    &\text{Finally, we truncate:} & W &= [\tilde U]_{R_V-\rank(C^T)}.
\end{align}
This construction is a generalization of the construction in \cite[Appendix C]{sanderse2020non}.

By inserting the velocity approximation \eqref{eq:velo expansion} into the momentum equation \eqref{eq:mom eq} and Galerkin-projecting, we find the ODE
\begin{align}
    \phi^T\Omega_h\ddt \phi \vv a(t) = \ddt \vv a(t)= \phi^T\vv{F_h^{CD}}(\phi \vv a(t),t) - \phi^T G_h \vv p_h(t).
    \label{eq:ins rom ode}
\end{align}
The right-hand side of this ODE contains a pressure term. We will eliminate this pressure term in the next section as a by-product of enforcing mass conservation.



\subsection{Mass conservation}
\label{sec:mass cons for arb velo basis}

Inserting the velocity approximation \eqref{eq:velo expansion} into the mass equation \eqref{eq:mass eq}, we find
\begin{align}
    M_h\phi \vv a(t) = 0.
    \label{eq:rom mass eq}
\end{align}
%
%
%
To satisfy mass conservation, we observe that the ROM \eqref{eq:velo expansion} satisfies the mass equation \eqref{eq:rom mass eq} for all $\vv a(t)\in\mbR^{R_V}$, if the ROM basis $\phi$ is divergence-free, i.e.\ $M_h\phi = 0$. 

To exploit this insight if $\phi$ is not divergence-free, we decompose the basis into a divergence-free component $\phi_0$ and a basis orthogonal to the space of divergence-free vectors, $\phi_\bot$.
For this purpose,
we compute the QR decomposition of $(M_h\phi)^T$,
\begin{align}
    (M_h\phi)^T = [Q_1^M \;\; Q_2^M]
    \begin{bmatrix}
    R_1^M \\ 0
    \end{bmatrix} = Q_1^M R_1^M
\end{align}
and decompose
\begin{align}
    \vv a(t) = Q_1^M \vv{a_1}(t) + Q_2^M \vv{a_2}(t).
    \label{eq:a decomp}
\end{align}
Then, we find $\phi \vv a(t) = \phi_0 \vv{a_2}(t) + \phi_\bot \vv{a_1}(t)$ with the divergence-free basis $\phi_0 = \phi Q_2^M$ and the basis $\phi_\bot = \phi Q_1^M$ orthogonal to the divergence-free subspace.
Inserting the decomposition \eqref{eq:a decomp} in the mass equation \eqref{eq:rom mass eq}, we find
\begin{align}
    M_h\phi \vv a(t) = (R_1^M)^T(Q_1^M)^T[Q_1^M\vv{a_1}(t) + Q_2^M \vv{a_2}(t)] = (R_1^M)^T\vv{a_1}(t).
\end{align}
Since $(R_1^M)^T$ has full column rank, the velocity approximation \eqref{eq:velo expansion} satisfies the mass equation \eqref{eq:rom mass eq}, if and only if we set $\vv{a_1}(t) = 0$. 
As a result, the velocity approximation simplifies to $V_r(t) = \phi Q_2^M \vv{a_2}(t) = \phi_0 \vv{a_2}(t)$.
The remaining coefficients $\vv{a_2}(t)$ are computed via the $(Q_2^M)^T$-premultiplied ODE \eqref{eq:ins rom ode},
\begin{align}
    (Q_2^M)^T\ddt Q_2^M \vv{a_2}(t) = \ddt \vv{a_2}(t) = (Q_2^M)^T\phi^T \vv{F_h^{CD}}(\phi_0 \vv{a_2}(t),t).
    \label{eq:ins velo only ode}
\end{align}
The pressure term is omitted because
$(Q_2^M)^T\phi^TG_h = \phi_0^TG_h = (M_h\phi_0)^T = 0$.

To solve this ODE, we propose to employ an energy-conserving Runge-Kutta method, e.g., \cite{sanderse2013energy}.


\subsection{Kinetic energy conservation}
\label{sec:kinetiv energy}
As in \cite{sanderse2020non}, we define the kinetic energy of the ROM as
\begin{align}
    K_r(t) = \frac{1}{2} \vv{V_r}(t)^T\Omega_h \vv{V_r}(t) = \frac{1}{2} \|\vv{a_2}(t)\|_2^2.
\end{align}
Using the ODE 
\eqref{eq:ins velo only ode}, we calculate
\begin{align}
    \ddt K_r(t) &= \vv{a_2}(t)^T\ddt \vv{a_2}(t) = \vv{a_2}(t)^T \phi_0^T \vv{F_h^{CD}}(\phi_0 \vv{a_2}(t),t), 
\end{align}
where there pressure term is eliminated due to the mass equation \eqref{eq:rom mass eq}.

As shown in \cite{sanderse2020non}, the right-hand side evaluates in absence of body forces to $-\nu\|Q_r \vv{a_2}\|_2^2$ with some matrix $Q_r\in\mbR^{N_V\times R_V}$. Hence, the kinetic energy is constant in the inviscid limit and decreases otherwise. Consequently, the ROM is nonlinearly stable.

In contrast, for the ROM proposed in \cite{carlberg2018conservative}, the kinetic energy evolution is given by
\begin{align}
    \ddt K_r(t) = -\nu \|Q_r \vv{a_2}(t)\|_2^2 + \vv{a_2}(t)^T \vv{F_h^{CD,*}}(\phi \vv{a_2}(t),t),
\end{align}
provided that we use the approach described in Section \ref{sec:mass cons for arb velo basis} to enforce mass conservation.
Because of the perturbation term $\vv{F_h^{CD,*}}(\phi \vv{a_2}(t),t)$ corresponding to \eqref{eq:perturbation term}, it is unclear whether the global kinetic energy increases over time and 
and whether nonlinear stability can be inferred.


\section{CONCLUSION}
\label{sec:conclusion}

In this article, we have proposed a novel approach to develop reduced order models (ROMs) for finite volume discretizations that preserve the underlying conservation law over arbitrary subdomains. The key element of this approach is the modification of the ROM basis in consideration of the predefined subdomains. This modification incorporates the subdomain conservation constraints into the test basis and thereby enforces their satisfaction.
The main advantage of our proposed ROM is that it does not introduce perturbation terms and can be simply defined via the conventional Galerkin ROM ODE.

For the special case that conservation quantities are constant in time over the chosen subdomain, we have shown that the ROMs can be computed more efficiently by introducing a physically motivated constant offset.

Furthermore, we have applied the novel subdomain conservation approach to an energy-conserving finite volume discretization of the incompressible Navier-Stokes equations. To enforce the mass equation on the ROM level, we have proposed a method based on a QR decomposition (which can also be generalized to other algebraic constraints). From the mass conservation, we have inferred the kinetic energy evolution which implies nonlinear stability of the ROM. This makes it stand out from existing subdomain conservative ROMs that involve perturbation terms, for which such an argument probably does not hold.
\section{OUTLOOK}
\label{sec:outlook}

In this paper we have not addressed the choice of subdomains and their effect on the accuracy of the resulting ROMs. 
We believe that subdomain conservation could improve the generalization accuracy of ROMs. In particular, we see a potential in choosing a small number of subdomains each consisting of a single finite volume. Consequently, the ROM would be simulated in the classical finite volume sense over these finite volumes, while the remaining finite volumes are simulated in the Galerkin ROM fashion.
This way, we can decompose the simulation domain into areas where we have confidence in the reduced basis and areas where we foresee that a high-fidelity method is required. 

This interpretation can be seen as a generalization of heterogeneously refined grids in classical numerical methods. Such heterogeneous grids have highly refined grid cells where they are believed to be necessary and coarse grid cells everywhere else. The size of the grid cells is either computed apriori based on prior knowledge or adaptively based on a posteriori error estimates. In a similar fashion we could determine the conservation subdomains based on prior knowledge, or based on the analysis of snapshot data of the FOM simulation.

An example of informative prior knowledge is a parametrized external force that acts locally in the domain. 
In the area where this force acts, the solution likely changes significantly over variation of the parameters. Hence, snapshot data of a given set of parametrizations might not generalize well to other parametrizations. Therefore, we suggest to compute the ROM via finite volumes in this area.

As an idea to determine conservation subdomains based on snapshot data, we suggest to compute in which finite volumes the conservation quantity changes the most over time. 
We interpret the conservation laws over these finite volumes as the most important ones to enforce and hence suggest to choose these finite volumes as conservation subdomains.

Apart from these ideas to improve the global accuracy, enforcing conservation over single finite volumes could be used to improve the ROM accuracy locally. This approach could be useful for applications that are characterized by a particular interest in simulation results only in a small portion of the computation domain, e.g., flow around an airfoil.


 
\bibliographystyle{elsarticle-num} 

\bibliography{library/references}

\end{document}